\documentclass[12pt]{article}
\usepackage{amssymb}

\newcommand{\text}{\mbox}

\newtheorem{theorem}{Theorem}[section]
\newtheorem{proposition}{Proposition}[section]
\newtheorem{lemma}{Lemma}[section]
\newtheorem{corollary}{Corollary}[section]
\newtheorem{definition}{Definition}[section]

\begin{document}

\author{Frank Oertel\\Department of Physics and 
Mathematics\\Zurich University of Applied Sciences, Winterthur (ZHW)\\CH--8401 Winterthur}
\title{On normed products of operator ideals which contain $\frak{L}_2$ as a factor}
\maketitle



\begin{abstract}
\noindent We investigate quasi-Banach operator ideal products 
$({\frak{A}}\circ{\frak{B}},\mathbf{A\circ B})$ which contain 
$(\frak{L}_2, \mathbf{L}_2)$ as a factor. In particular, we ask 
for conditions which guarantee that $\mathbf{A\circ B}$ is even a norm 
if each factor of the product is a $1$-Banach ideal. In doing so, we reveal
the strong influence of the existence of such a norm in relation to 
the accessibility of the product ideal and the structure of its factors.
\newline
\newline
{\noindent}\textit{Key words and phrases:} Accessibility, Banach spaces,
cotype 2, finite rank operators, Hilbert space factorization, operator ideals, 
principle of local reflexivity
\newline
\newline
{\noindent}\textit{2000 AMS Mathematics Subject Classification:} primary 46A32, 
46M05, 47L20; secondary 46B07, 46B10, 46B28. 

\end{abstract}

\section{Introduction}

This paper is devoted to an investigation of the 
normability of operator ideal products which contain
$(\frak{L}_2, \mathbf{L}_2)$ as a factor (where $\frak{L}_2$ denotes the class of all 
operators which factor through a Hilbert space). 
It seems that $1$-Banach operator ideal products play a fundamental 
role in the search for maximal Banach ideals which do not satisfy a
transfer of the norm estimation in the classical principle of local
reflexivity to their ideal norm (cf. \cite{oe5}). This problem (which still is open) 
originated from the objective to facilitate the search for a non--accessible maximal 
normed Banach ideal (which is the same as a
non--accessible finitely generated tensor norm in the sense of Grothendieck)
(cf. \cite{oe1}). Later, in 1993, Pisier
constructed a counterexample (cf. \cite[31.6.]{df}). Since each
right--accessible maximal Banach ideal $(\frak{A},\mathbf{A})$ even
satisfies such a \textit{principle of local reflexivity for operator ideals}, 
Pisier`s counterexample of a non--accessible maximal Banach ideal 
naturally lead to the search
for counterexamples of maximal Banach ideals
$({\frak{A}}_0,{\mathbf{A}}_0)$ for which the conjugate
$({\frak{A}}_0^\Delta ,{\mathbf{A}}_0^\Delta )$ is not left-accessible, 
implying surprising relations between the
existence of a norm on product operator ideals of type 
$\frak{B}\circ \frak{L}_2$, the extension of finite rank operators with respect to a
suitable operator ideal norm and the principle of local reflexivity for
operator ideals (cf. \cite{oe5}). The basic objects, connecting these
different aspects, are product operator ideals with the property (I) and
the property (S), introduced by Jarchow and Ott (see \cite{jo}). In the widest
sense, a product operator ideal $\frak{A}\circ \frak{B}$ has the property
(I), if
\[
(\frak{A}\circ \frak{B})\cap \frak{F}=(\frak{A}\cap \frak{F})\circ \frak{B}
\]
and the property (S), if
\[
(\frak{A}\circ \frak{B})\cap \frak{F}=\frak{A}\circ (\frak{B}\cap \frak{F})
\]
(where $\frak{F}$ denotes the class of all finite rank operators) so that
each finite rank operator in $\frak{A}\circ \frak{B}$ is the composition of
two operators, one of which is of finite rank. Since each operator ideal
which contains $\frak{L}_2$ as a factor, has both, the property (I) and the
property (S), Hilbert space factorization crystallized out as a fundamental
key in these investigations.

\section{The framework}

In this section, we introduce the basic notation and terminology which we
will use throughout in this paper. We only deal with Banach spaces and most
of our notations and definitions concerning Banach spaces and operator
ideals are standard. We refer the reader to the monographs \cite{df}, \cite
{djt} and \cite{p1} for the necessary background in operator ideal theory
and the related terminology. Infinite dimensional Banach spaces over the
field $\Bbb{K\in \{R},\Bbb{C\}}$ are denoted throughout by $W,X,Y$ and $Z$
in contrast to the letters $E,F$ and $G$ which are used for finite
dimensional Banach spaces only. Denote for given Banach spaces $X$ and $Y$
\[
\text{FIN}(X):=\{E\subseteq X\mid E\in \text{FIN}\}\hspace{0.2cm}\text{and}%
\hspace{0.2cm}\text{COFIN}(X):=\{L\subseteq X\mid X/L\in \text{FIN}\}\text{,}
\]
where FIN stands for the class of all finite dimensional Banach spaces. 
The space of all operators (continuous
linear maps) from $X$ to $Y$ is denoted by $\frak{L}$(X,Y), and for the
identity operator on $X$, we write $Id_X$. The collection of all finite rank
(resp. approximable) operators from $X$ to $Y$ is denoted by $\frak{F}$(X,Y)
(resp. $\overline{\frak{F}}$(X,Y)), and $\frak{E}$(X,Y) indicates the
collection of all operators, acting between finite dimensional Banach spaces
$X$ and $Y$ (elementary operators). The dual of a Banach space $X$ is
denoted by $X^{\prime }$, and $X^{\prime \prime }$ denotes its bidual $%
(X^{\prime })^{\prime }$. If $T\in \frak{L}$(X,Y) is an operator, we
indicate that it is a metric injection by writing $T:X\stackrel{1}{%
\hookrightarrow }Y$, and if it is a metric surjection, we write $T:X%
\stackrel{1}{\twoheadrightarrow }Y$. If $X$ is a Banach space, $E$ a finite
dimensional subspace of $X$ and $K$ a finite codimensional subspace of $X$,
then $B_X:=\{x\in X:\Vert x\Vert \leq 1\}$ denotes the closed unit ball, $%
J_E^X$ $:E\stackrel{1}{\hookrightarrow }X$ the canonical metric injection
and $Q_K^X:X\stackrel{1}{\twoheadrightarrow }X\diagup K$ the canonical
metric surjection. Finally, $T^{\prime }\in {\frak{L}}(Y^{\prime },X^{\prime
}) $ denotes the dual operator of $T\in \frak{L}$(X,Y).

If $({\frak{A}},\mathbf{A})$ and $({\frak{B}},\mathbf{B})$ are given
quasi--Banach ideals, we will use throughout the shorter notation $({\frak{A}}%
^{\text{d}},\mathbf{A}^{\text{d}})$ for the dual ideal and the abbreviation ${\frak{A}}\stackrel{%
1}{=}{\frak{B}}$ for the isometric equality $(\frak{A},\mathbf{A})=(\frak{B},%
\mathbf{B})$. We write $\frak{A}\subseteq \frak{B}$ if, regardless of the
Banach spaces $X$ and $Y$, we have ${\frak{A}}(X,Y)\subseteq {\frak{B}}(X,Y)$.
The metric inclusion $(\frak{A},\mathbf{A})\subseteq (\frak{B}, \mathbf{B})$ 
is often shortened by $\frak{A}\stackrel{1}{\subseteq }\frak{B}$. 

Given quasi--Banach ideals $(\frak{A},\mathbf{A)}$ and
$(\frak{B},\mathbf{B)} $, let $({\frak{A}}\circ
{\frak{B}},\mathbf{A\circ B})$ be the corresponding product ideal
and $({\frak{A}}\circ {\frak{B}}^{-1},{\mathbf{A}}\circ {\mathbf{B}}^{-1})$
(resp. $({\frak{A}}^{-1}\circ
{\frak{B}},{\mathbf{A}}^{-1}\circ{\mathbf{B}}))$ the
corresponding ''right--quotient'' (resp. ''left--quotient''). 
Important examples are $(\frak{A}^{\text{min}}, \mathbf{A}^{\text{min}}):=
({\overline{\frak{F}}} \circ {\frak{A}} \circ {\overline{\frak{F}}}, 
\Vert \cdot \Vert \circ {\mathbf{A}} \circ \Vert \cdot \Vert)$ 
(the \textit{minimal kernel} of $(\frak{A},\mathbf{A)}$) and 
$(\frak{A}^{\text{max}}, \mathbf{A}^{\text{max}}):=
({\overline\frak{F}}^{-1} \circ {\frak{A}} \circ {\overline\frak{F}}^{-1}, 
\Vert \cdot \Vert^{-1} \circ {\mathbf{A}} \circ \Vert \cdot \Vert^{-1})$ 
(the \textit{maximal hull} of $(\frak{A},\mathbf{A)}$). 
$({\frak{A}}^{\text{inj}},{\mathbf{A}}^{\text{inj}})$ denotes 
the \textit{injective hull}\emph{\ }of $\frak{A}$, the unique smallest 
injective quasi--Banach ideal which contains $({\frak{A}},\mathbf{A})$, 
and $({\frak{A}}^{\text{sur}},{\mathbf{A}}^{\text{sur}})$, 
the \textit{surjective hull} of $\frak{A}$, is the unique smallest 
surjective quasi--Banach ideal which contains $({\frak{A}},\mathbf{A})$. 

In addition to the maximal Banach ideal
$({\frak{L}},\mathcal{\Vert \cdot \Vert })$ we mainly will be concerned with the
maximal Banach ideals $({\frak{I}},\mathbf{I})$ (integral operators),
$({\frak{L}}_2,{\mathbf{L}}_2)$ (Hilbertian
operators), $({\frak{D}}_2,{\mathbf{D}}_2)\stackrel{1}{=}
({\frak{L}}_2^{\displaystyle{*}},{\mathbf{%
L}}_2^{\displaystyle{*}})\stackrel{1}{=}{\frak{P}}_2^{\text{d}}\circ {\frak{P}}_2
$\/ ($2$--dominated operators), $({\frak{P}}_p,{\mathbf{P}}_p)
$ (absolutely $p$--summing operators),
$1\leq p\leq \infty ,\frac
1p+\frac 1q=1$, $({\frak{L}}_\infty ,{\mathbf{L}}_\infty )\stackrel{1}{=}
({\frak{P}}_1^{\displaystyle{*}},{\mathbf{P}}_1^{\displaystyle{*}})$ and
$({\frak{L}}_1,{\mathbf{L}}_1)%
\stackrel{1}{=}({\frak{P}}_1^{\displaystyle{*}\text{d}},
{\mathbf{P}}_1^{\displaystyle{*}\text{d}})$.\\

Since we will use them throughout in this paper, let us recall the important
notions of the conjugate operator ideal (cf. \cite{glr}, \cite{jo} and \cite
{oe2}) and the adjoint operator ideal (all details can be found in the
standard references \cite{df} and \cite{p1}). Let $(\frak{A},\mathbf{A})$ be
a quasi--Banach ideal.

\begin{itemize}
\item  Let ${\frak{A}}^\Delta (X,Y)$ be the set of all $T\in {\frak{L}}(X,Y)$
which satisfy
\[
{\mathbf{A}}^\Delta (T):=\sup \{\mid \text{tr}(TL)\mid \text{ }\mid L\in {\frak{F}}%
(Y,X),{\mathbf{A}}(L)\leq 1\}<\infty .
\]
Then a Banach ideal $({\frak{A}}^\Delta ,{\mathbf{A}}^\Delta )$ is
obtained (here, $\text{tr}(\cdot )$ denotes the usual trace for finite
rank operators). It is called the \textit{conjugate ideal} of
$(\frak{A},\mathbf{A})$.

\item  Let ${\frak{A}}^{\displaystyle{*}}(X,Y)$ be the set of all $T\in {\frak{L}}(X,Y)$
which satisfy
\[
{\mathbf{A}}^{\displaystyle{*}}(T):=\sup \{\mid \text{tr}(TJ_E^XSQ_K^Y\}\mid \text{
}\mid E\in \text{FIN}(X),K\in \text{COFIN}(Y),{\mathbf{A}}(S)\leq 1\}<\infty .
\]
Then a Banach ideal $({\frak{A}}^{\displaystyle{*}},{\mathbf{A}}^{\displaystyle{*}})$ is
obtained. It is called the \textit{adjoint operator ideal }of
$(\frak{A},\mathbf{A})$.
\end{itemize}

A deeper investigation of relations between the Banach ideals $({\frak{A}}^\Delta ,%
{\mathbf{A}}^\Delta )$ and $({\frak{A}}^{\displaystyle{*}},{\mathbf{A}}^{\displaystyle{*}})$
needs the help of an important local property, known as
accessibility, which can be viewed as a local version of
injectivity and surjectivity. All necesary details about
accessibility of operator ideals and its applications can be
found in \cite {df}, \cite{oe2}, \cite{oe3} and \cite{oe4}. So
let us recall :

\begin{itemize}
\item  A quasi--Banach ideal $(\frak{A},\mathbf{A})$ is called \textit{%
right--accessible}, if for all $(E,Y)\in $ FIN $\times $ BAN, operators $%
T\in {\frak{L}}(E,Y)$ and $\varepsilon >0$ there are $F\in $ FIN$(Y)$ and $%
S\in {\frak{L}}(E,F)$ so that $T=J_F^YS$ and ${\mathbf{A}}(S)\leq
(1+\varepsilon ){\mathbf{A}}(T)$.

\item  $(\frak{A},\mathbf{A})$ is called \textit{left--accessible}, if for
all $(X,F)\in $ BAN $\times $ FIN, operators $T\in {\frak{L}}(X,F)$ and $%
\varepsilon >0$ there are $L\in $ COFIN$(X)$ and $S\in
{\frak{L}}(X/L,F)$ so that $T=SQ_L^X$ and ${\mathbf{A}}(S)\leq
(1+\varepsilon ){\mathbf{A}}(T)$.

\item  A left--accessible and right--accessible quasi--Banach ideal is called 
\textit{accessible}.

\item  $(\frak{A},\mathbf{A})$ is \textit{totally accessible}, if for every
finite rank operator $T\in {\frak{F}}(X,Y)$ acting between Banach spaces $X$%
, $Y$ and $\varepsilon >0$ there are $(L,F)\in $ COFIN$(X)\times
$ FIN$(Y)$ and $S\in {\frak{L}}(X/L,F)$ so that $T=J_F^YSQ_L^X$
and ${\mathbf{A}}(S)\leq (1+\varepsilon ){\mathbf{A}}(T)$.
\end{itemize}

Let us recall the following important results on accessibility (for a 
detailed proof cf. \cite{df}, \cite{oe1} and \cite{oe4}):

\begin{theorem} Let $(\frak{A},\mathbf{A})$ be a Banach
ideal. Then $(\frak{A}^{\text{min}}, \mathbf{A}^{\text{min}})$ always is accessible and 
$(\frak{A}^{{\displaystyle*}\Delta},\mathbf{A}^{{\displaystyle*}\Delta})$
is right-accessible. If in addition $(\frak{A}, \mathbf{A})$ is maximal, then 
$(\frak{A}, \mathbf{A})$ is right-accessible if and only 
if $(\frak{A}^{\displaystyle*},\mathbf{A}^{\displaystyle*})$ is left-accessible.
\end{theorem}

\begin{theorem}
Let $(\frak{A},\mathbf{A})$ be a maximal Banach ideal.
\begin{itemize}
\item[$(i)$] $(\frak{A},\mathbf{A})$ is right-accessible if and 
only if $\frak{A}^{\displaystyle*}\circ\frak{A} \stackrel{1}{\subseteq }\frak{I}$.
\item[$(ii)$] $(\frak{A},\mathbf{A})$ is totally accessible if and only 
if $\frak{A}^{\displaystyle*}\stackrel{1}{=}\frak{A}^{\Delta}$.
\end{itemize}
\end{theorem}

Pisier's counterexample $({\frak{A}}_P,{\mathbf{A}}_P)$ 
shows the existence of maximal Banach
ideals which neither are left nor right--accessible (cf. \cite{df}, 31.6). 
However, accessibility of a quasi--Banach ideal at least can be 
transmitted to its regular hull:

\begin{proposition}
Let $(\frak{A},\mathbf{A)}$ be an arbitrary quasi--Banach ideal. 
If $(\frak{A},\mathbf{A)}$ is right--accessible (resp. totally--accessible), 
then the regular hull $(\frak{A}^{\text{reg}},\mathbf{A^{\text{reg}}})$ 
is also right--accessible (resp. totally--accessible).
\end{proposition}

\textsc{Proof: }Let $\epsilon >0$, $X$, $Y$ be Banach spaces and $T\in \frak{%
F}(X,Y)$ an arbitrary finite rank operator. Assume that $\frak{A}$ is
totally accessible or that $X\in FIN$ and $\frak{A}$ is right--accessible.
In both cases, there exists a finite dimensional Banach space $F\in
FIN(Y^{\prime \prime })$ and an operator $S\in \frak{L}(X,F)$, so that $%
j_YT=J_F^{Y^{\prime \prime }}S$ and 
\[
\mathbf{A}(S)<(1+\epsilon )\cdot \mathbf{A}(j_YT)=(1+\epsilon )\cdot \mathbf{%
A}^{\text{reg}}(T)\text{.} 
\]
Due to the classical principle of local reflexivity for linear operators,
there exists an operator $W\in \frak{L}(F,Y)$ so that $\mathbf{\Vert }W%
\mathbf{\Vert }<1+\epsilon $ and $j_YWz=J_F^{Y^{\prime \prime }}z$ for all $%
z\in F$ which satisfy $J_F^{Y^{\prime \prime }}z\in j_Y(Y)$. Let $x\in X$
and put $z:=Sx$. Then $J_F^{Y^{\prime \prime }}z=$ $j_YTx\in j_Y(Y)$, which
therefore implies that $j_YWSx=J_F^{Y^{\prime \prime }}z=$ $j_YTx$. Now,
factor $W$ canonically through a finite dimensional subspace $G$ of $Y$ so
that $W=J_G^YU$ and $\mathbf{\Vert }U\mathbf{\Vert }<1+\epsilon $.
Consequently, $T=WS=J_G^Y(US)$, and 
\[
\mathbf{A}^{\text{reg}}(US)<(1+\epsilon )^2\cdot \mathbf{A}^{\text{reg}}(T)\text{.} 
\]
Hence, $\frak{A}^{\text{reg}}$ is right--accessible (in each of the both cases). In
the case of $\frak{A}$ being totally accessible, the operator $S$ even can
be chosen as $S=S_0Q_K^X$, where $K\in COFIN(X)$ and $S_0\in \frak{L}%
(X\diagup K,F)$ so that 
\[
\mathbf{A(}S_0)<(1+\epsilon )\cdot \mathbf{A}^{\text{reg}}(T)\text{,} 
\]
and the proof is finished.$\blacksquare $

\section{Normed operator ideal products}

Let $(\frak{A},\mathbf{A)}$ be a $p$-normed Banach ideal $(0<p \leq 1)$
and $(\frak{B},\mathbf{B)}$ be a $q$-normed Banach ideal $(0<q \leq 1)$. Then
the product $({\frak{A}}\circ{\frak{B}},\mathbf{A\circ B})$ is a $r$-normed Banach
ideal, where $1/r:=1/p +1/q$ (see \cite{p1}, 7.1.2). Even if $p=1$ and $q=1$,
$({\frak{A}}\circ{\frak{B}},\mathbf{A\circ B})$ in general is a $1/2$-Banach ideal only;
$\mathbf{A\circ B}$ need not to be a norm. However, if one of the operator ideals 
is closed (such as e.g., $(\overline{\frak{F}}, \Vert \cdot \Vert)$, $({\frak{K}},\Vert \cdot \Vert)$ or 
$({\frak{W}},\Vert \cdot \Vert$)), then we may formulate a positive result (cf. \cite{b}):

\begin{proposition}
Let $(\frak{A},\mathbf{A)}$ and $(\frak{B},\mathbf{B)}$ be two quasi-Banach ideals. Let 
$0 <p \leq1$. Then, in each of the following cases, 
$({\frak{A}}\circ{\frak{B}},\mathbf{A\circ B})$ is a $p$-Banach ideal:
\begin{itemize}
\item[$(i)$] $(\frak{A},\mathbf{A})$ is a $p$-Banach ideal and $(\frak{B},\mathbf{B}) = 
(\frak{B},\Vert \cdot \Vert)$ is closed,
\item[$(ii)$] $(\frak{A},\mathbf{A}) = (\frak{A},\Vert \cdot \Vert)$ is closed 
and $(\frak{B},\mathbf{B})$ is a $p$-Banach ideal.
\end{itemize}
\end{proposition}

\textsc{Proof:} It is sufficient to prove the case (i); (ii) follows similarly. So let
$(\frak{A},\mathbf{A})$ be $p$-normed and $(\frak{B},\Vert \cdot \Vert)$ be 
closed. Let $X$ and $Y$ be arbitrary Banach spaces and $T_1$, $T_2 \in 
({\frak{A}}\circ{\frak{B}})(X,Y)$. It remains to show that
\begin{center} 
$(\mathbf{A} \circ \Vert \cdot \Vert)^p(T_1 +T_2) \leq 
(\mathbf{A} \circ \Vert \cdot \Vert)^p(T_1) + 
(\mathbf{A} \circ \Vert \cdot \Vert)^p(T_2)$.
\end{center}
Let $\epsilon >0$. Then there exist Banach spaces $Z_1, Z_2$ and operators
$R_1 \in \frak{A}(Z_1, Y)$,  $R_2 \in \frak{A}(Z_2, Y)$, $S_1 \in \frak{B}(X, Z_1)$,
$S_2 \in \frak{B}(X, Z_2)$ so that $T_1=R_1 S_1$, $T_2=R_2 S_2$, $\Vert S_1 \Vert \leq 1$,
$\Vert S_2 \Vert \leq 1$, $\mathbf{A}(R_1) \leq
(1+\epsilon) \cdot (\mathbf{A} \circ \Vert \cdot \Vert)(T_1)$ and 
$\mathbf{A}(R_2) \leq
(1+\epsilon) \cdot (\mathbf{A} \circ \Vert \cdot \Vert)(T_2)$. 
We now consider the Banach space $W:=l_\infty (Z_1, Z_2)$ consisting of 
all elements $(z_1,z_2)\in
Z_1\times Z_2$ so that $\Vert (z_1,z_2)\Vert _\infty :=\max (\Vert z_1\Vert,
\Vert z_2\Vert )<\infty $. Let $J_i : Z_i \stackrel{1}{\hookrightarrow }W$ be the canonical
injections and $Q_i : W\stackrel{1}{\twoheadrightarrow }Z_i$ the corresponding 
canonical surjections ($i=1,2$). Then $S:=J_1S_1+J_2S_2 \in \frak{B}(X, W)$ and 
$\Vert Sx \Vert_\infty = \Vert(S_1x, S_2x)\Vert_\infty = \max (\Vert S_1x\Vert,
\Vert S_2x\Vert ) \leq \Vert x \Vert$ for all $x \in X$. Hence, 
$\Vert S \Vert \leq 1$. Put $R:=R_1Q_1 +R_2Q_2$. Then $R \in \frak{A}(W, Y)$ and
$(\mathbf{A}(R))^p \leq (\mathbf{A}(R_1))^p + (\mathbf{A}(R_2))^p$. The construction 
therefore implies that $T_1+T_2 =RS$ and
$(\mathbf{A} \circ \Vert \cdot \Vert)^p (T_1 +T_2) =(\mathbf{A} \circ \Vert \cdot \Vert)^p (RS)
\leq (\mathbf{A}(R) \cdot \Vert S \Vert)^p \leq (\mathbf{A}(R))^p  \leq (\mathbf{A}(R_1))^p 
+ (\mathbf{A}(R_2))^p \leq (1+\epsilon)^p \cdot ((\mathbf{A} \circ \Vert \cdot \Vert)^p (T_1) +
(\mathbf{A} \circ \Vert \cdot \Vert)^p (T_2))$,
and the proof is finished. $\blacksquare$\\

An immediate (non-trivial) consequence is the

\begin{corollary}
Let $(\frak{A},\mathbf{A)}$ be a $p$-Banach ideal $(0 < p \leq 1)$. Then 
$(\frak{A}^{\text{min}},\mathbf{A}^{\text{min}})$ also is a $p$-Banach ideal.
\end{corollary}

Unfortunately, we still cannot present explicite sufficient 
criteria which show the existence of (an equivalent) ideal norm on product ideals in the
general case.
It seems to be much more easier to show that a certain product ideal cannot be a
normed one by using arguments which involve trace ideals and the ideal of
nuclear operators (the smallest Banach ideal). Even more holds: if $\frak{A}%
\circ {\frak{L}}_2$ is a 1--Banach ideal for certain operator ideals $\frak{A}$%
, then $\frak{A}\circ {\frak{L}}_2$ \textit{is not} right--accessible (cf. theorem 3.4)! 
To investigate more carefully the general case, we recall an important 
factorization property for finite rank operators which 
had been introduced by Jarchow and Ott in their paper \cite{jo}. 
It not only turns out to be very useful for an investigation of local structures in 
(product) operator ideals; this factorization property was used as 
the main tool to show that $\frak{L}_\infty $ and 
$\frak{L}_1$ are \textit{not} totally accessible -- answering an open question of Defant
and Floret (see \cite{df}, 21.12 and \cite{oe5}). So, let us recall the definition of this
factorization property and its implications:

\begin{definition}[Jarchow/Ott]
Let $(\frak{A},\mathbf{A)}$ and $(\frak{B},\mathbf{B)}$ be arbitrary
quasi--Banach ideals. Let $L\in \frak{F}(X,Y)$ an arbitrary finite rank
operator between two Banach spaces $X$ and $Y$. Given $\epsilon >0$, we can
find a Banach space $Z$ and operators $A\in \frak{A}(Z,Y)$, $B\in \frak{B}%
(X,Z)$ so that $L=AB$ and 
\[
\mathbf{A}(A)\cdot \mathbf{B}(B)\leq (1+\epsilon )\cdot (\mathbf{A\circ B})(L)%
\text{.} 
\]

\begin{enumerate}
\item[$(i)$]  If the operator $A$ is of finite rank, we say that $\frak{A}%
\circ \frak{B}$ has the property\textit{\ }(I).

\item[$(ii)$]  If the operator $B$ is of finite rank, we say that $\frak{A}%
\circ \frak{B}$ has the property\textit{\ }(S).
\end{enumerate}
\end{definition}

Important examples are the following (see \cite{jo}, lemma 2.4.):

\begin{itemize}
\item  If $\frak{B}$ is injective, or if $\frak{A}$ contains $\frak{L}_2$ as
a factor, then $\frak{A}\circ \frak{B}$ has the property\textit{\ }(I).

\item  If $\frak{A}$ is surjective, or if $\frak{B}$ contains $\frak{L}_2$
as a factor, then $\frak{A}\circ \frak{B}$ has the property\textit{\ }(S).
\end{itemize}

Since $\frak{L}_2\circ \frak{A}$ is injective for every quasi--Banach ideal $%
(\frak{A},\mathbf{A})$ (see \cite{oe4}, lemma 5.1.), $\frak{B}\circ \frak{L}%
_2\circ \frak{A}$ therefore has the property (I) as well as the property
(S), for all quasi--Banach ideals $(\frak{A},\mathbf{A})$ and $(\frak{B},\mathbf{B})$. 
Such ideals are exactly those which contain $\frak{L}_2$ as factor -- in the sense 
of \cite{jo}.

\begin{theorem}
Let $(\frak{A},\mathbf{A})$ be a maximal Banach ideal. Then both,
the maximal $\frac 12$--Banach ideal ${\frak{A}}^{\text{inj}}\circ
{\frak{L}}_2$ and the injective hull of the maximal $\frac
12$--Banach ideal $\frak{A}\circ {\frak{A}}_2$ are totally
accessible.
\end{theorem}

\textsc{Proof:} Since every Hilbert space has the metric
approximation property and since
${\frak{A}}^{\text{inj}}\stackrel{1}{=}({\frak{A}}^{\text{inj}})^{\displaystyle{**}}$ is
right--accessible, an easy calculation shows that
\begin{equation}
{\frak{A}}^{\text{inj}}\circ
{\frak{L}}_2\stackrel{1}{=}({\frak{A}}^{\text{inj}})^{{\displaystyle*}\Delta }\circ
{\frak{L}}_2\text{.}
\end{equation}
Since $({\frak{A}}^{\text{inj}})^{{\displaystyle*}\Delta }$ is right--accessible, the
total
accessibility of ${\frak{L}}_2$ and the property (S) of the product ideal $(%
{\frak{A}}^{\text{inj}})^{{\displaystyle*}\Delta }\circ {\frak{L}}_2$ even imply that $({\frak{A}}%
^{\text{inj}})^{{\displaystyle*}\Delta }\circ {\frak{L}}_2$ is totally accessible (cf.
\cite[proposition 4.1]{oe5}). Hence, ${\frak{A}}^{\text{inj}}\circ
{\frak{L}}_2$ is totally accessible (due to $(1)$), and in
particular we obtain that $({\frak{A}}\circ
{\frak{L}}_2)^{\text{inj}}\stackrel{1}{=}({\frak{A}}^{\text{inj}}\circ
{\frak{L}}_2)^{\text{inj}} $ is totally accessible.$\blacksquare$\\

Now, let $(\frak{A},\mathbf{A})$ be a maximal Banach ideal so that
${\mathbf{L}}_2\circ \mathbf{A}$ is even a norm on the (maximal) product ideal $({\frak{L}}%
_2\circ {\frak{A}},{\mathbf{L}}_2\circ {\mathbf{A}})$. Then ${\frak{A}}^{\displaystyle{*}}\stackrel{%
1}{\subseteq }({\frak{L}}_2\circ {\frak{A}})^{\displaystyle{*}}\stackrel{1}{\subseteq }{\frak{L}}%
_\infty $ (cf. \cite[proposition 5.1.]{oe4}) and ${\frak{L}}_\infty \stackrel{1%
}{=}{\frak{P}}_1^\Delta \stackrel{1}{\subseteq }{\frak{N}}^\Delta
$. Given Banach spaces $X$ and $Y$ so that both, $X^{\prime }$
and $Y$ have cotype 2, \cite[theorem 4.9.]{pi} tells us, that any
finite rank operator $L\in {\frak{F}}(Y,X)$ satisfies
\[
{\mathbf{N}}(L)\leq (2{\mathbf{C}}_2(X^{\prime })\cdot
{\mathbf{C}}_2(Y))^{\frac 32}\cdot {\mathbf{D}}_2(L)\text{.}
\]
Hence,
\[
{\frak{N}}^\Delta (X,Y)\subseteq {\frak{D}}_2^\Delta (X,Y)\stackrel{1}{=}{\frak{L}}%
_2(X,Y)\text{,}
\]
and we have proven a rather surprising fact (revealing the strong influence
of a \textit{norm} on an operator ideal product):

\begin{theorem}
Let $(\frak{A},\mathbf{A})$ be a maximal Banach ideal so that the
product ideal $({\frak{L}}_2\circ {\frak{A}},{\mathbf{L}}_2\circ
\mathbf{A})$ is normed. Let $X$ and $Y$ be arbitrary Banach
spaces so that both, $X^{\prime }$ and $Y $ have cotype 2. Then
\[
{\frak{A}}^{\displaystyle{*}}(X,Y)\stackrel{1}{\subseteq }({\frak{L}}_2\circ {\frak{A}}%
)^{\displaystyle{*}}(X,Y)\subseteq {\frak{L}}_2(X,Y)
\]
and
\[
{\mathbf{L}}_2(T)\leq (2{\mathbf{C}}_2(X^{\prime })\cdot
{\mathbf{C}}_2(Y))^{\frac 32}\cdot ({\mathbf{L}}_2\circ
{\mathbf{A}})^{\displaystyle{*}}(T)\leq (2{\mathbf{C}}_2(X^{\prime })\cdot
{\mathbf{C}}_2(Y))^{\frac 32}\cdot {\mathbf{A}}^{\displaystyle{*}}(T)
\]
for all operators $T\in {\frak{A}}^{\displaystyle{*}}(X,Y)$.
\end{theorem}

To maintain the previous statement, even a permutation of the factors $\frak{%
A}$ and ${\frak{L}}_2$ in the product ${{\frak{L}}}_2\circ
\frak{A}$ is allowed:

\begin{theorem}
Let $(\frak{A},\mathbf{A})$ be a maximal Banach ideal so that the
product ideal $({\frak{A}}\circ {\frak{L}}_2,{\mathbf{A\circ
L}}_2)$ is normed. Let $X$ and $Y$ be arbitrary Banach spaces so
that both, $X^{\prime }$ and $Y$ have cotype 2. Then
\[
{\frak{A}}^{\displaystyle{*}}(X,Y)\stackrel{1}{\subseteq }({\frak{A}}\circ {\frak{L}}%
_2)^{\displaystyle{*}}(X,Y)\subseteq {\frak{L}}_2(X,Y)
\]
and
\[
{\mathbf{L}}_2(T)\leq (2{\mathbf{C}}_2(X^{\prime })\cdot
{\mathbf{C}}_2(Y))^{\frac 32}\cdot ({\mathbf{A\circ
L}}_2{\mathbf{)}}^{\displaystyle{*}}(T)\leq (2{\mathbf{C}}_2(X^{\prime })\cdot
{\mathbf{C}}_2(Y))^{\frac 32}\cdot {\mathbf{A}}^{\displaystyle{*}}(T)
\]
for all operators $T\in {\frak{A}}^{\displaystyle{*}}(X,Y)$.
\end{theorem}

\textsc{Proof:} Let $(\frak{A},\mathbf{A})$ and $X$, $Y$ be as
before and let ${\frak{A}}\circ {\frak{L}}_2$ be normed. Then
${\frak{A}}\stackrel{1}{=}{\frak{A}}^{\text{dd}}$, and ${\frak{A}}\circ
{\frak{L}}_2$ is a maximal (and therefore a
regular) Banach ideal (cf. \cite[lemma 4.3]{oe5}). Since the injective $%
\frac 12$--Banach ideal ${\frak{L}}_2\circ {\frak{A}}^{\text{d}}$ is also
regular (cf. \cite[lemma 5.1]{oe4}), an easy calculation shows
that
\[
({\frak{A}}\circ {\frak{L}}_2)^{\text{d}}\stackrel{1}{=}{\frak{L}}_2\circ
{\frak{A}}^{\text{d}}
\]
and\footnote{%
In particular, it follows that ${\frak{A}}\circ {\frak{L}}_2$ is
surjective (cf. \cite[8.5.9.]{p1}).}
\[
{\frak{A}}\circ {\frak{L}}_2\stackrel{1}{=}({\frak{L}}_2\circ
{\frak{A}}^{\text{d}})^{\text{d}}\text{.}
\]
Since ${\mathbf{A\circ L}}_2$ is a norm, $({\mathbf{A\circ L}}%
_2)^{\text{d}}$ obviously is a norm too. Hence, if we apply the previous theorem to
the normed product ideal $({\frak{A}}\circ {\frak{L}}_2)^{\text{d}}\stackrel{1}{=}{\frak{L}}%
_2\circ {\frak{A}}^{\text{d}}$, we obtain
\[
{\frak{A}}^{\displaystyle{*}\text{d}}(X,Y)\stackrel{1}{\subseteq }({\frak{L}}_2\circ {\frak{A}}%
^{\text{d}})^{\displaystyle{*}}(X,Y)\stackrel{1}{=}({\frak{A}}\circ
{\frak{L}}_2)^{\displaystyle{*}\text{d}}(X,Y)\subseteq {\frak{L}}_2(X,Y)\text{,}
\]
and
\[
{\mathbf{L}}_2(T)\leq C\cdot ({\mathbf{A\circ
L}}_2)^{\displaystyle{*}}(T^{\prime })\leq C\cdot {\mathbf{A}}^{\displaystyle{*}}(T^{\prime })
\]
for \textit{all} operators $T\in {\frak{A}}^{\displaystyle{*}\text{d}}(X,Y)$ (where $C:=(2{\mathbf{C}}%
_2(X^{\prime })\cdot {\mathbf{C}}_2(Y))^{\frac 32}$). Now, since
$Y$ has the same coptype as its bidual $(Y^{\prime })^{\prime }$
with identical cotype
constants (cf. \cite[corollary 11.9]{djt}), the proof is finished.$%
\blacksquare$\\

Let $(\frak{A},\mathbf{A})$ be a given ultrastable quasi--Banach
ideal so that $({\frak{A}}\circ {\frak{L}}_2,{\mathbf{A\circ
L}}_2)$ is right--accessible.
Our aim is to show that in this case $({\frak{A}}\circ {\frak{L}}_2,{\mathbf{%
A\circ L}}_2)$ and $({\frak{L}}_2\circ {\frak{A}}^{\displaystyle{*}},{\mathbf{L}}_2\circ {\mathbf{A}}%
^{\displaystyle{*}})$ \textit{both together} cannot be normed.\textit{\ }To this end, we
need a lemma which is of its own interest:

\begin{lemma}
Let $({\frak{A}}_0,{\mathbf{A}}_0)$ be a maximal Banach ideal so that space$(%
{\frak{A}}_0)$ contains a Banach space without the approximation
property. Then there does not exist a maximal Banach ideal
$(\frak{C},\mathbf{C})$ so
that ${\frak{C}}\circ {\frak{L}}_\infty $ has the property (I) and ${\frak{C}}%
\subseteq {\frak{A}}_0^{-1}\circ {\frak{P}}_1$.
\end{lemma}

\textsc{Proof:} Assume that the statement is false. Then there
exists a (maximal) Banach ideal $(\frak{A},\mathbf{A})$ so that
${{\frak{A}}}_0\subseteq {{\frak{P}}}_1\circ
({\frak{A}}^{\displaystyle{*}})^{-1}\stackrel{1}{=}({\frak{A}}^{{\displaystyle*}\Delta })^{\text{inj}}
$. Due to the assumed property (I) of ${\frak{A}}^{\displaystyle{*}}\circ
{\frak{L}}_\infty $, the proof of theorem 3.4 in \cite{oe5} shows
that even $(({\frak{A}}^{{\displaystyle*}\Delta
})^{\text{inj}})^{\text{dd}}\stackrel{1}{\subseteq }({\frak{A}}^{\text{inj}})^{{\displaystyle*}\Delta }\stackrel{1%
}{\subseteq }{\frak{N}}^\Delta $. Since ${\frak{A}}_0$ was assumed
to be a
maximal Banach ideal, we therefore obtain ${\frak{A}}_0\stackrel{1}{=}{\frak{A}}%
_0^{\text{dd}}\stackrel{1}{\subseteq }{\frak{N}}^\Delta $ which is a contradiction.$%
\blacksquare $

\begin{corollary}
Let $({\frak{A}}_0,{\mathbf{A}}_0)$ be a maximal Banach ideal so that space$(%
{\frak{A}}_0)$ contains a Banach space without the approximation property. If $%
({\frak{A}}_0^{-1}\circ {\frak{P}}_1)\circ {\frak{L}}_\infty $ has
the property (I), ${\frak{A}}_0$ is not left--accessible.
\end{corollary}

\begin{theorem}
Let $(\frak{B},\mathbf{B})$ be an ultrastable quasi--Banach ideal so that $%
{\frak{B}}\subseteq {\frak{L}}_\infty $. If ${\frak{B}}\circ
{\frak{L}}_2$ is right--accessible, ${\frak{B}}\circ {\frak{L}}_2$
cannot be a $1-$ Banach ideal.
\end{theorem}

\textsc{Proof:} Assume that the statement is false and put ${\frak{B}}_0:=(%
{\frak{L}}_\infty \circ {\frak{L}}_2)^{\displaystyle{*}}$ and ${\frak{A}}:=({\frak{B}}\circ {\frak{L}}%
_2)^{\displaystyle{*}}$. Then
\[
{\frak{A}}^{\displaystyle{*}}\stackrel{1}{=}({\frak{B}}\circ {\frak{L}}_2)^{\text{max}}\stackrel{1}{=}(%
{\frak{B}}\circ
{\frak{L}}_2)^{\text{reg}}\stackrel{1}{=}{\frak{B}}^{\text{reg}}\circ
{\frak{L}}_2
\]
is right--accessible (due to proposition 2.1) and contains ${\frak{L}}%
_2$ as a factor so that in particular ${\frak{A}}^{\displaystyle{*}}\circ
{\frak{L}}_\infty $ has the property (I). Since
${\frak{B}}\subseteq {\frak{L}}_\infty $,
\[
{\frak{B}}_0\circ {\frak{A}}^{\displaystyle{*}}\subseteq {\frak{A}}\circ {\frak{A}}^{\displaystyle{*}}\stackrel{1}{%
\subseteq }{\frak{I}}\stackrel{1}{\subseteq }{\frak{P}}_1\text{,}
\]
and it follows that ${\frak{A}}^{\displaystyle{*}}\subseteq
{\frak{B}}_0^{-1}\circ {\frak{P}}_1$. Since $Id_P\in {\frak{B}}_0$
(cf. \cite[proposition 4.4]{oe5}), lemma 3.1 leads to a
contradiction.$\blacksquare$\\

Now let us assume that $(\frak{B},\mathbf{B})$ is even a \textit{maximal
Banach ideal} so that ${\frak{B}}\subseteq {\frak{L}}_\infty $. If ${\frak{B}}%
\circ {\frak{L}}_2$ were normed, then ${\frak{B}}\circ
{\frak{L}}_2$ would be a
maximal \textit{and }surjective Banach ideal, implying that ${\frak{P}}_1^{\text{d}}%
\stackrel{1}{=}{\frak{I}}^{\text{sur}}\stackrel{1}{=}({\frak{N}}^{\text{max}})^{\text{sur}}%
\stackrel{1}{\subseteq }{\frak{B}}\circ {\frak{L}}_2$. Hence,
\begin{equation}
({\frak{B}}\circ {\frak{L}}_2)^{\displaystyle{*}}\stackrel{1}{\subseteq }{\frak{P}}_1^{\text{d}\displaystyle{*}}%
\stackrel{1}{=}{\frak{L}}_1\text{.}
\end{equation}
Since ${\frak{B}}\subseteq {\frak{L}}_\infty $, it even follows that ${\frak{P}}%
_1\stackrel{1}{=}{\frak{L}}_\infty ^{\displaystyle{*}}\subseteq ({\frak{B}}\circ {\frak{L}}_2)^{\displaystyle{*}}%
\stackrel{1}{\subseteq }{\frak{L}}_1$ which is a contradiction
(cf. \cite[27.2.]{df}). So, in this case we obtain a stronger
result:

\begin{theorem}
Let $(\frak{B},\mathbf{B})$ be a maximal Banach ideal so that ${\frak{B}}%
\subseteq {\frak{L}}_\infty $. Then ${\frak{B}}\circ {\frak{L}}_2$
cannot be a $1-$ Banach ideal.
\end{theorem}

\begin{corollary}
Let $(\frak{A},\mathbf{A})$ be a quasi-Banach ideal. If $(\frak{A},\mathbf{A})$ is 
a maximal Banach ideal or if $(\frak{A},\mathbf{A})$ is regular and 
$({\frak{A}}\circ {\frak{L}}_2,{\mathbf{A\circ L}}_2)$ is right-accessible, then
$({\frak{A}}\circ {\frak{L}}_2,{\mathbf{A\circ L}}_2)$ and $({\frak{L}}_2\circ {\frak{A}}^{\displaystyle{*}},{\mathbf{L}}_2
\circ {\mathbf{A}}^{\displaystyle{*}})$ both together cannot be normed.
\end{corollary}

\textsc{Proof:} Let $({\frak{A}}\circ {\frak{L}}_2,{\mathbf{A\circ L}}_2)$ be a $1$-Banach ideal.
Due to our assumption, ${\frak{A}} \not\subseteq {\frak{L}}_\infty $. If the injective quasi-Banach ideal 
$({\frak{L}}_2\circ {\frak{A}}^{\displaystyle{*}},{\mathbf{L}}_2\circ {\mathbf{A}}^{\displaystyle{*}})$ were also a normed one, then
we would obtain ${\frak{P_1}} \stackrel{1}{\subseteq} 
{\frak{L}}_2\circ {\frak{A}}^{\displaystyle{*}} \stackrel{1}{\subseteq} {\frak{A}}^{\displaystyle{*}}$ 
(cf. \cite{oe4}, proposition 5.1) -- a contradiction.$\blacksquare$\\

\noindent \textbf{Acknowledgements:} The author would like to thank
the Organizing and Program Committee of 
\textit{Functional Analysis 2000 Valencia} for the invitation, the perfect
establishment of this magnificent event and the friendly
atmosphere. In particular, he would like to thank G. Arango, K.
Floret, H. Junek, A. Mart\'{i}nez--Abej\'{o}n, E. A. S\'{a}nchez
P\'{e}rez, A. Sofi, V. Tarieladze, O. J. Tylli and R. Vidal for a
fruitful exchange of ideas and stimulating discussions. In addition, 
he would like to thank an anonymous referee for some very useful 
remarks.


\end{document}